%


\input amstex      
\documentstyle{amsppt}      
      
\magnification=1200      
\pagewidth{6.5truein}      
\pageheight{9.0truein}

\define\C{\Bbb C}      
\define\CC{{\Cal C}}    
\define\R{\Bbb R}

\define\eps{\epsilon}

\define\ls{\mathrel{\mathop{<}\limits_{\sim}}}      
\define\gs{\mathrel{\mathop{>}\limits_{\sim}}}

\define\CS{{\Cal S}} 
\define\CR{{\Cal R}}
      
\redefine\Re{\operatorname{Re}}      
\redefine\Im{\operatorname{Im}}

\define\aut{\operatorname{Aut}}

\topmatter

\title      
On a Domain in $\C^2$ With generic piecewise smooth Levi-flat boundary and
Non-compact Automorphism Group
\endtitle      
\rightheadtext{Non-Compact Automorphism Groups}      
\author      
Siqi Fu and Bun Wong
\endauthor      
      
      
\address      
Department of Mathematics, \ \ \ Texas A \& M University, \ \ \
College Station, \ \ TX 77840
\endaddress

\email       
sfu\@ math.tamu.edu 
\endemail      

\address
Department of Mathematics,\ \ \ University of California, \ \ \
Riverside,\ \ CA\ 92521
\endaddress      
      
\email
wong\@ math.ucr.edu
\endemail

      
      
\endtopmatter

\document      
    
\subhead \S 1. Introduction \endsubhead    
    
In this paper, we continue our investigation of domains with non-compact
automorphism groups in [FW].  We shall prove the following:

\proclaim{Main Theorem} If $D$ is a simply-connected domain in 
$\C^2$ with generic piecewise smooth Levi-flat boundary and non-compact
automorphism group, then $D$ is biholomorphic to the bidisc.
\endproclaim

The boundary $bD$ of a  bounded domain $D$  in $\C^n$ is called {\it
piecewise smooth} if there exists a neighborhood
$U$ of $\overline{D}$ and $\rho_k \in \CC^\infty (U)$,
$1\le k \le m$, such that $D =\{ q\in U; \
\rho_k(q)<0, 1\le k\le m \}$ and $d\rho_{k_1}\wedge \ldots
\wedge d\rho_{k_l} \not= 0$
on $\cap_{j=1}^{l}S_{k_j}$ for any distinct 
$k_1, \ldots, k_l \in \{ 1, \ldots, m \}$, 
where $S_j=\{q\in U; \  \rho_j(q)=0 \}$.
It is called {\it generic piecewise smooth} if 
$\partial\rho_{k_1}\wedge \ldots \wedge \partial\rho_{k_l} \not= 0$
on $\cap_{j=1}^{l}S_{k_j}$.
The boundary $bD$  is called {\it (generic) piecewise smooth Levi-flat} if
each $S_j$ is  in additional Levi-flat (See Section 3).  We will call
$\{\rho_j; \ \  1\le j \le m\}$ a {\it defining system} of
$D$ and each $S_j$ a {\it defining hypersurface} of $D$.

When $D$ is convex, the above result was obtained by 
K.-T. Kim [Kim1] (see [W2] for related results). Kim's 
proof  uses a refine
version of the rescaling method introduced by Frankel [Fra].
 It was proved
by Pinchuk [P] that a homogeneous bounded domain with piecewise
smooth boundary is biholomorphic to a product of balls.  Note that
the non-compact condition in the above theorem is  weaker than the 
homogeneous condition in Pinchuk's result.  In the latter case, 
one can choose a special boundary accumulation point that has
properties similar to those possessed by 
a strictly pseudoconvex boundary point.  
See [Kod1, 2] and [CS] for results along this line. 
We remark that 
the simply-connected condition on $D$ cannot be dropped.  
For example,
the product of a disc and an annulus has generic piecewise smooth Levi-flat
boundary and non-compact automorphism group. 
However, it is not biholomorphic to the bidisc.  For motivation
and background on the subject, we refer readers to 
[W1,2], [R], [GK1-3], [BP1,2], [Kim1,2],
[FIK],  and  references therein. 
     
\bigskip  

\noindent{\bf Acknowledgment:}  The first author thanks Professors
H. Boas, S. Krantz and E. Straube for their supports.  He is especially
thankful to Professor E. Straube for stimulating conversations on the
subject.

\bigskip 

\subhead \S 2. Preliminaries \endsubhead

Let $D$ be a bounded domain in $\C^n$.
Let $\Delta$ be the unit disc and $\Delta^n$   
the unit n-polydisc.     
Let $H(D_1, D_2)$ be the family of holomorphic   
mappings from $D_1$ to $D_2$. 
Let $T(D)$ be the holomorphic tangent bundle of $D$.  
we will identify $T(D)$ with $D\times \C^n$.    
   
The Kobayashi-Royden metric 
$ F^K_D\colon T(D)$ $\to \R^+\cup \{0\}$    
is defined by     
$$    
F^K_D(z, v)=\inf \left\{ 1/\lambda;  
\quad \text{there exits $f\in H(\Delta, D)$    
with $f(0)=z$, $f'(0)=\lambda v$},\ \lambda > 0\  \right\}    
$$    
for  $(z, v)\in T(D)$. 

For $f(z)=(f_1(z_1, \ldots, z_n), \ldots, f_m(z_1, \ldots, z_n))$, 
we denote by $f'(z)$ the $m\times n$ Jacobian matrix $\left( \partial 
f_j/\partial z_k\right)$.  
The Eisenman-Kobayashi measure on $D$ is defined by
$$    
M^E_D(z) =\inf\left\{ \frac1{|\det f'(0)|^2}; \    
f\in H(\Delta^n, D),\ \  f(0)=z\right\} .
$$    
The Carath\'{e}odory measure $M_D^C$ on $D$ is 
defined by
$$    
M_D^C(z)=\sup \{ |\det f'(z) |^2; 
\ f\in H(D, \Delta^n),\ \ f(z)=0 \}.    
$$    

We list some well-known properties of the Kobayashi-Royden
metric and the invariant measures.  

\proclaim{Lemma 2.1} Let $D$, $D_1$ and  $D_2$ be bounded domains in 
$\C^n$.

\item{(1)} 
If $f\in H( D_1, D_2)$, then 
$F^K_{D_1}(p, v)\ge F^K_{D_2}(f(p), f_{*p}(v))$   
for $(p, v)\in T(D_1)$. 

\item{(2)} 
If $f\in H(D_1, D_2)$, then    
$M_{D_1}(z) \ge M_{D_2}(f(z))|\det f'(z)|^2$,     
where $M_D$ is either of the invariant measures.

\item{(3)} Let $\pi\colon D_1\to D_2$ be a covering mapping, then 
$F^K_{D_1}(z, v)=F^K_{D_2}(\pi(z), \pi_*(v))$ and
$M^E_{D_1}(z)$ = $M^E_{D_2} (\pi (z)) |\det \pi'(z) |^2$.    

\item{(4)} $M^E_D(z)\ge M_D^C(z)$ for all $z\in D$. 
If $M^E_D(z_0)= M_D^C(z_0)$ for some $z_0\in D$, then    
$D$ is biholomorphic to $\Delta^n$.

\endproclaim    

For more information on
invariant metrics and measures, we refer readers to
[Kob], [Kr], [GW1,2], and [JP].

Let $\Gamma_\theta=\{ re^{i\alpha}; \ \ r>0, \ 
\pi-\theta < \alpha < \pi +\theta \}$ 
be the cone with vertex at   
the origin and angles between  its edges and the  
negative $\Re z$-axis $\theta$.  Let 
$\Delta(a, r)$ be the disc with center $a$ and radius
$r$.  Let $\Delta_\eps =\Delta(0, \eps)$ and
$\Gamma_\theta^\eps =\Gamma_\theta \cap \Delta_\eps$.   

\proclaim{Lemma 2.2} Let $\theta_1$ and $\theta_2$ be two numbers   
such that $0<\theta_1<\theta_2<\pi$. 
Let $z_j =r_j e^{i(\pi +\alpha_j)}$, $-\theta_1
<\alpha_j < \theta_1 $,  be a  
sequence in $\Gamma_{\theta_1}$. 
Suppose that $z_j\to 0$ and $\alpha_j\to \alpha<\theta_1$.  
Then for any positive numbers $\eps_1$ and  $\eps_2$, 
   
$$   
\lim_{j\to \infty} \frac{M_{\Gamma_{\theta_1}^{\eps_1}}(z_j)}   
{M_{\Gamma_{\theta_2}^{\eps_2}}(z_j)}    
=   
\left(\frac{\theta_2\cos \frac{\alpha\pi}{2\theta_2}}
{\theta_1\cos\frac{\alpha\pi} 
{2\theta_1}}\right)^2   
$$   
where $M$ is either of the invariant measures.   
\endproclaim   
     
\demo{Proof}  Note that in this case the two invariant measures are
identical.  In fact, they are  the square of the Poincar\'{e}  metric 
in the unit direction.  Using the conformal mapping 
$$
z\mapsto \left( \frac{1 + i\left(-\frac{z}{\eps}\right)^{\pi/2\theta}}
{1 - i\left(-\frac{z}{\eps}\right)^{\pi/2\theta}}\right)^2
$$
that maps $\Gamma^\eps_\theta$ onto the upper-half plane and the explicit
formula for Poincar\'{e} metric of the upper-half plane, one obtains
by straight-forward but tedious calculations that for $z=re^{i(\pi + \phi)}$
with $-\theta < \phi < \theta$, 
$$
M_{\Gamma^\eps_\theta}(z)
=\frac{1}{\left(\theta\cos\frac{\pi\phi}{2\theta}\right)^2}\cdot
\frac{\pi^2\left( 1 + (\frac{r}{\eps})^{\frac{\pi}{\theta}} 
-2(\frac{r}{\eps})^{\frac{\pi}{2\theta}}\sin  \frac{\pi\phi}{2\theta}\right)
\left( 1 + (\frac{r}{\eps})^{\frac{\pi}{\theta}} 
+2(\frac{r}{\eps})^{\frac{\pi}{2\theta}}\sin  \frac{\pi\phi}{2\theta}\right)
}{4r\left( 1- (\frac{r}{\eps})^{\frac{\pi}{\theta}}\right)^2} .
$$
The lemma then follows from the above formula. \qed
\enddemo

\subhead \S 3 Levi-flat hypersurfaces \endsubhead

Let $S$ be a smooth hypersurface in $\C^n$.  Let $\rho$ be a defining
function for $S$, {\it i.e.}, there exists a neighborhood $U$ of $S$
such that $\rho\in C^\infty (U)$, $S=\{z\in U; \ \ \rho(z)=0 \}$, and
$d\rho (z)\not= 0 $ on $S$. The hypersurface 
$S$ is called {\it pseudoconvex} if the Levi-form
$$
\sum_{j, k=1}^n \frac{\partial^2 \rho}{\partial z_j\partial \bar z_k} (z) 
X_j\bar X_k \ge 0
$$
for all $z\in S$ and $X\in \C^n$ with $\sum_{j=1}^n \frac{\partial\rho}
{\partial z_j}(z) X_j = 0$.  It is called {\it Levi-flat} if
the ``$\ge$'' sign in the above inequality is replaced by the ``$=$'' sign.   

It is well-known that a Levi-flat hypersurface  is locally
foliated by complex manifolds of codimension 1 (cf. [Fre]).  In particular,
if $S$ is a Levi-flat hypersurface in $\C^2$ 
and $p\in S$, then there exists a neighborhood
$U$ of $p$ and a diffeomorphism $g_t(\zeta)=g(t, \zeta)\colon (-1, 1)\times
\Delta \to  S\cap U$ such that $g_t(\zeta)$ is holomorphic in $\zeta$.  Each
open Riemann surface $g_t(\Delta)$ is a leaf of the foliation.  The following
lemma shows that if $S$ contains an affine disc, 
then one can piece together the local
foliations to obtain a foliation of $S$ in a neighborhood of any smaller
disc of the given disc.
Since we cannot find a reference for such a result, we provide 
the details of the proof, which is inspired in part by the work of Barrett and
Forn{\ae}ss [BaF].

\proclaim{Lemma 3.1} Let $M=\Delta\times\{ 0 \}$ 
and let $S$ be a Levi-flat hypersurface
containing $M$.  Then for any $\delta$ with $0<\delta<1$
there exists a neighborhood $N_\delta$ of
$\Delta_\delta\times\{ 0\}$ and a diffeomorphism $\Phi(t, \zeta)$
from $(-1, 1)\times \Delta_\delta$ onto 
$S\cap N_\delta$ such that 
$\Phi (t, \zeta)=(\zeta, \varphi(t, \zeta))$ where $\varphi(t, \zeta)$
is holomorphic in $\zeta$ and $\varphi (0, \zeta)=0$ for $\zeta\in
\Delta_\delta$.
\endproclaim

\demo{Proof} It follows from Theorem 8 in [DF] (see Lemma 3.3 below) 
that after a change
of coordinate system in a neighborhood of $M$, the real normal of 
$S$ is constant on $M$. Assume that
the normal direction is the positive real $\Re w$-direction. Then
$S$ is given in a neighborhood of $\Delta_{\delta}\times \{ 0\}$
by
$$
\Re w + r(z, \bar z, \Im w) = 0  \tag 2.1
$$
where $|r(z, \bar z, \Im w)|\le C|\Im w|^2$ for some constant $C>0$.  

We now cover the closure of $\Delta_{\delta}$ by finitely many
discs $\Delta(z_j, r_j)$, $0\le j \le m$,  so that
$S$ has a foliation by open Riemann surfaces in a
neighborhood $U_j$ of $\Delta(z_j, 2r_j)\times \{0\}$.
 
Without loss of generality, we assume that $\Delta(z_0, r_0)$
contains the origin. Let $g(t, \zeta)=(g^1(t, \zeta), g^2(t, \zeta))
\colon (-1, 1)\times\Delta\to S\cap U_0$ be a foliation of $S\cap U_0$.  
After reparametrization, 
we may assume that $g(0, 0)=0$. Since $g(0, \Delta)\in S$,
it follows from (2.1) that $|\Re g^2(0, \zeta)|\le C |\Im g^2(0, \zeta)|^2$.
Thus, $g^2(0, \zeta)\equiv 0$ and $g^1(0, \zeta)$ is one-to-one for $\zeta
\in \Delta$. Let $a\in (0, 1)$ be sufficiently closed to 1 such that
$\Delta(z_0, r_0)\subset\subset g^1(0, \Delta_a)$.  It is easy to
see from Rouch\'{e}'s theorem that $g^1_t(\zeta)=g^1(t, \zeta)$
is one-to-one on $\Delta_a$ for sufficiently small $t$.   
 
Let $\Phi_0(t, \zeta)=(\zeta,
\varphi_0 (t, \zeta))\colon (-\eps_0, \eps_0)\times \Delta(z_0, r_0)\to
S$ where $\varphi_0(t, \zeta)=g^2_t \circ
(g^1_t)^{-1}(\zeta)$ and $\eps_0$ is a
sufficiently small positive constant.  After a reparametrization in $t$, we 
may assume that $\Phi_0(t, 0)=(0, -r(0, 0, t)+ it)$. 
It is clear that $\Phi_0 (t, \zeta)$
gives a foliation of $S$ in a neighborhood of $\Delta(z_0, r_0)\times\{ 0\}$.

We now show how to extend $\Phi_0(t, \zeta)$ to 
obtain a foliation of $S$ in a neighborhood of $\Delta_\delta\times \{0\}$.
Suppose that $\Delta(z_0, r_0)\cap \Delta(z_1, r_1)\not=\emptyset$.  
Let $p\in \Delta (z_0, r_0)\cap \Delta (z_1, r_1)$.
Let $\Phi_1(t, \zeta)=(\zeta, \varphi_1(t, \zeta))\colon  
(-\eps_1, \eps_1)\times \Delta(z_1, r_1)\to S$ be defined as
in the previous paragraph such that
$\Phi_1(t, p)=(p, \varphi_0(t, p))$.  Here $\eps_1$ is a
sufficiently small constant. 
It follows from the uniqueness of 
the foliation that $\Phi_0(t, \zeta)=\Phi_1(t, \zeta)$ for
$t\in (-\eps_0, \eps_0)\cap (-\eps_1, \eps_1)$ and $\zeta\in 
\Delta(z_0, r_0)\cap \Delta(z_1, r_1)$.  
Piecing together the foliations of $S$ in a neighborhood
of $\Delta(z_j, r_j)\times \{ 0\}$ in this manner, 
we obtain a smooth foliation 
of $S$ in a neighborhood of $\Delta_\delta$ in a form of
$$
\Phi (t, \zeta)=(\zeta, \varphi(t, \zeta))\colon  
(-\eps', \eps')\times \Delta_\delta\to S
$$
where $\Phi(t, 0)=(0, -r(0, 0, t) +i t)$ and $\eps' > 0$ is a 
sufficiently small constant.
After a reparametrization in $t$, we may assume that $\eps'=1$.
\qed

\enddemo

\proclaim{Lemma 3.2} Let $D$ be a bounded domain in $\C^2$ with
piecewise smooth Levi-flat boundary.  Let $S$ be one of the
defining hypersurfaces
of $D$.   Let $\hat D$ be a subdomain of $D$ and let $g\colon 
\hat D \to bD$ be a holomorphic map.  
If $g(\hat D)\cap S\not=\emptyset$, then $g(\hat D)\subseteq S$. 
\endproclaim

\demo{Proof} We assume that $g$ is a non-constant map.
Let $M=g(\hat D)\cap S$.  Then $M$ is a closed subset
of $g(\hat D)$.  We now prove that it is also an open subset.  
Let $p\in M$ and $p=g(q)$ for some $q\in \hat D$.   
After a change of coordinate system, we may assume that $p=(0, 0)$,
the leaf of $S$ through $p$ is locally parameterized by $\zeta\mapsto
 (\zeta, 0)$,  and $S$ is locally defined by 
$$
\rho(z, w)=\Re w + \Im w \cdot r(z, \bar z, \Im w) 
$$
where $r(z, \bar z, \Im w)=O(|z|+|\Im w|)$.  Write $g=(g^1, g^2)$.  It follows
from $\rho(g^1, g^2)\le 0$ that 
$\Re g^2(z, w) \le | \Im g^2(z, w) |$ for $(z, w)$ sufficiently
closed to $q$.  Since $g^2(q)=0$, it follows from the open mapping property 
of non-constant holomorphic functions that $g^2$ is identically zero near
$q$.  Since $M$ contains the set $\{ (g^1(z, w), 0); \ \text{for
all $(z, w)$ near $q$} \}$, it is also an open subset of $g(\hat D)$.
Therefore $g(\hat D)=M\subseteq S$.   
\qed
\enddemo

Since the above-mentioned  result of Diederich and Forn{\ae}ss 
[DF, Theorem~8] plays an important role in this paper, we
reformulate it here for reader's convenience.  

\proclaim{Theorem 3.3 (Diederich-Forn{\ae}ss)}  Let $D$ be a bounded domain
in $\C^2$ such that $M=\Delta\times \{ 0 \}\subseteq bD$. Assume that
$bD$ is smooth and pseudoconvex in a neighborhood of $M$ with a local
defining function $r=r(z, w)$. Assume further that the outward
normal direction of $bD$ at the origin is the positive $\Re w$-axis. 
Let $v(z)=\arg \frac{\partial r}{\partial w} (z, 0)$. 
Then $v(z)$ defines a harmonic function on $\Delta$.  Furthermore, if 
$u(z)$ is a harmonic conjugate of $v(z)$ and $h(z)=\exp(-u(z)
+iv(z))$, then the outward normal of $bD$ is constant on $M$
in the new coordinates $(z', w')$ defined by $z'=z$ and $w'=wh(z)$.     
\endproclaim
\bigskip

\subhead \S 4. Proof of the Main Theorem, Part I \endsubhead

For the rest of the paper, we will use $D$ to denote
a simply-connected domain in $\C^2$ with generic
piecewise smooth Levi-flat boundary and
non-compact automorphism group $\aut(D)$.
Let $\{g_k\}\subseteq \aut(D)$.  Suppose
that $g_k(q)\to bD$ as $k\to \infty$ for some $q\in D$.  
After passing to a subsequence,
we may assume that $g_k \to g$ local uniformly on $D$. 
It follows from Cartan's theorem (cf. [Na1, pp. 78]) that
$g\colon D\to bD$.

Let $\CS =bD \cap (\cup_{j\not= k} (S_j\cap S_k))$ be the set
of singular boundary points. Let $\CR=bD\setminus \CS$ be the set
of regular boundary points.  

In this section, we prove the main theorem for the case when 
$g(D)\cap \CR \not= \emptyset$.

Let $p=g(q)\in g(D)\cap \CR$.  Assume that $p\in S$ where $S$
is one of the defining hypersurfaces of $D$ with defining function $\rho$.
Then there exists a neighborhood
$U$ of $p$ so that $D\cap U =\{ (z, w)\in U;\ 
\rho (z, w) < 0 \}$.
It follows from Lemma 3.2 that $g(D)\subseteq 
S$.

\proclaim{Lemma 4.1} With above notations and conditions, 
$g\not\equiv$ constant. 
Furthermore, if 
$\hat D$ is a relatively compact subdomain of
$D$, then $g(\hat D)$ is a locally closed open Riemann 
surface.
\endproclaim

\demo{Proof}  Let $q_j = g_j (q)$
and let $p_j$ be the projection of $q_j$ onto the boundary
in the direction of the outward normal of $bD$ at $p$. 
Let $M_j$ be the leaf on $S$ that passes through $p_j$ and
let $v_j$ be the unit complex tangent vector of $M_j$ at $p_j$.  It is easy to
see that  
$$
\|(g_j^{-1})'(q_j)v_j\|\ls F^K_D (q, (g^{-1}_j)'(q_j)v_j)
= F^K_D (q_j, v_j) \ls 1 
$$
for all $j$. Let $\tilde v_j =
(g_j^{-1})'(q_j)v_j/ \|(g_j^{-1})'(q_j)v_j\|$.  Passing to a subsequence,
we may assume that $\tilde v_j $ converges to a unit vector $\tilde v$.
Since $|g_j'(q)\tilde v_j | \gs 1$, we have that $|g'(q) \tilde v|
\gs 1$.  Therefore $g\not\equiv$ constant.  

Let $q'\in \hat D$ and $p'=g(q')$.  After a change of local coordinate
near $p'$, we assume that $p'=(0, 0)$, the positive $\Re w$-axis is
the outward direction of $bD$ at $p'$, and the leaf of $S$ through $p'$
is locally parameterized by $\zeta\mapsto (\zeta, 0)$.  Let $g=(g^1, g^2)$.
It follows from the proof of Lemma 3.2 that $g^2$ is identically zero 
in a neighborhood of $q'$.  Therefore $g^1\not\equiv$ constant. Thus 
$g(\hat D)\cap (\Delta_\eps\times\Delta_\eps) 
\supseteq \Delta_\eps \times \{0 \}$ for 
sufficiently small $\eps >0$.   
It remains to prove that after possible shrinking 
of $\eps$,  $g(\hat D)\cap (\Delta_\eps\times \Delta_\eps)
= \Delta_\eps\times \{ 0\}$.  Suppose that
this is not true. Then there exists
a sequence $p'_j\in g(\hat D)$ such that $p'_j\to p'$ and
the second coordinate of $p'_j$ is not zero. Let $q'_j\in \hat D$
be such that $p'_j=g(q'_j)$.  After passing to a subsequence,
we assume that $q'_j\to \tilde q\in D$.  Since $g(\tilde q)=p'$, we have
again that $g^2(z, w)\equiv 0$ for $(z, w)$ in a neighborhood of $\tilde q$.
This contradicts to the assumption.   \qed   

\enddemo

\proclaim{Lemma 4.2} Let $\Omega$ be a bounded domain such that
$M=\Delta\times \{ 0 \} \subseteq b\Omega$ and
$S=b\Omega\cap (\Delta\times \Delta_{\eps_0})$ is smooth and Levi-flat for some
$\eps_0 >0$. Let $\delta\in (0, 1)$, $\eps >0$, 
and $U_{\delta\eps}=\Delta_\delta\times \Delta_\eps$.
Then for any sequence $\{ q_j\}$ in $\Omega$ that tends to $(0, 0)$,
$$ 
\varlimsup_{\eps\to 0}\varlimsup_{j\to \infty}   
\frac{ M^C_{\Omega\cap U_{\delta\epsilon}}(q_j)}    
{ M^E_{\Omega\cap U_{\delta\epsilon}}(q_j)}  =1  .     
$$    
\endproclaim

\demo{Proof} By Theorem 3.3, we can choose a
coordinate system in a neighborhood of $M$ such that
the outward normal of $b\Omega$ is the positive $\Re w$-axis
for points on $M$. Then $S$ is given in a neighborhood of
$\Delta_\delta\times \{ 0 \}$ by a defining function
of form
$$
\rho (z, w)=\Re w + r(z, \bar z, \Im w)
$$
where $|r(z, \bar z, \Im w)|\le C |\Im w|^2$ for some constant
$C>0$.
 
Let $N_{\delta}$  be a neighborhood
of $\Delta_{\delta}\times\{ 0\}$ and $\Phi(t, \zeta)=(\zeta,
\varphi(t, \zeta))\colon (-1, 1)\times \Delta_{\delta}\to
S\cap N_{\delta}$ be the diffeomorphism constructed in Lemma 3.1.

Write $q_j=(z_j, a_j + ib_j)$.  Let $t_j\in
(-1, 1)$ satisfy $\varphi(t_j, z_j)=-r(z_j, \bar z_j, b_j)
+ ib_j$.  Such  $t_j$ is uniquely determined and 
$t_j\to 0$.  It is clear that $L_j=\Phi(t_j, \Delta_\delta)$ is the leaf of
through the projection
of $q_j$ onto $S$ in the direction of the positive $\Re w$--axis.

Let $v_j(z)=\arg \frac{\partial \rho}{\partial w}(z, \varphi(t_j, z))$
where the argument takes the principal branch.  Since $\frac{\partial
\rho}{\partial w}(z, 0)= 1/2$ for $|z|<\delta$, $v_j(z)$ is well-defined
for sufficiently large $j$.  It follow from Theorem~3.3 that $v_j(z)$
is harmonic.  Let $u_j(z)$ be its harmonic conjugate such that $u_j(0)=0$.
Let $h_j(z')=\exp(-u_j(z')
+iv_j(z'))$.  Let $F_j\colon (z, w)\mapsto (z', w')$ be defined
by
$z'=z$ and $w'=(w-\varphi(t_j, z))h_j(z)$. Then
$$
F_j(\Omega\cap U_{\delta\eps})=\left\{ (z', w'); \ \
|z'|<\delta, \ |w'/h_j(z') +\varphi(t_j, z')|<\eps, \ \tilde \rho (z', w') <0 
\right\}
$$
where $\tilde\rho(z', w')=\rho(F^{-1}_j (z', w'))$.  Since
$$
\tilde\rho(z', 0)=0, \frac{\partial \tilde\rho}{\partial z'}(z', 0)=0, 
\ \text{and } 
\frac{\partial \tilde\rho}{\partial w'}(z', 0)=
|\frac{\partial \rho}{\partial w}(z', \varphi(t_j, z'))|e^{u_j(z')} >0
$$
for $|z'|<\delta$, therefore,
in the $(z', w')$-coordinates, $L_j=\Delta_\delta\times
\{ 0 \}$ and the outward normal of $b\Omega$ is the positive $\Re w'$-axis
for points on $L_j$.  On the other hand, there exists a constant
$C>0$ such that $|\varphi(t_j, z)|\le C|t_j|$ and 
$1/C\le |h_j(z)| \le C$ for $|z|<\delta$ and sufficiently large $j$.
Therefore, there exist constants $a>0, b>1$ such that
$$
\gather
\Delta_\delta\times\{|w'|<\eps/b;\  \Re w' + a|\Im w'|^2 < 0\} \\
\subseteq F_j (\Omega\cap U_{\delta\eps}) \subseteq \\
\Delta_\delta\times\{|w'|< b\eps;\  \Re w' +a|\Im w'|^2 < 0\}
\endgather
$$
for sufficiently small $\eps$ and sufficiently large $j$.  Thus for any two 
angles
${\theta_1}\in (0, \pi/2)$ and ${\theta_2}\in (\pi/2, \pi)$, we have
$\Delta_\delta\times \Gamma^{\frac{\eps}{b}}_{\theta_1}\subseteq
F_j(\Omega\cap U_{\delta\eps}) \subseteq \Delta_\delta\times
\Gamma^{b\eps}_{\theta_2}$ provided $\eps$ is sufficiently small.  Let $q'_j=
F_j(q_j)$ and $w'_j$ be the second coordinate of $q'_j$.  
It follows that
$$
\align
1\ge \frac{   M^C_{\Omega\cap U_{\delta\epsilon}}(q_j) }    
{  M^E_{\Omega\cap U_{\delta\epsilon}}(q_j) }   
&\ge
\frac{   M^C_{\Delta_\delta\times \Gamma^{b\eps}_{{\theta_2}}}(q'_j) }    
{  M^E_{\Delta_\delta\times \Gamma^{\frac{\eps}{b}}_{{\theta_1}}}(q'_j) } \\    
&= \frac{ M^C_{\Gamma^{b\eps}_{{\theta_2}}}(w'_j) }
{ M^E_{\Gamma^{\frac{\eps}{b}}_{{\theta_1}}}(w'_j) }  . \\
\endalign
$$
By Lemma 2.2, the last term can be chosen to be as 
close to 1 as we wish provided
$j\to \infty$, $\eps\to 0$,  ${\theta_1}\to (\pi/2)^-$, and 
${\theta_2}\to (\pi/2)^+$.  \qed
\enddemo

The following lemma is well-known and its proof is elementary.

\proclaim{Lemma 4.3} Let $D$ be a simply-connected domain. Let
$D_n$ be subdomains of $D$ such that $D_j\subset\subset D_{j+1}$ and
$\cup_{j=1}^n D_j =D$.  Then $D_j$ is also simply-connected 
for sufficiently large $j$.
\endproclaim

We are now in position to prove the main theorem in the case
when $g(D)\cap\CR\not= \emptyset$.  The proof uses ideas
from our previous work [FW].  We shall keep the notation
and setup as at the beginning of this section.  

Let $D_1$ and $D_2$ be simply-connected subdomains of $D$
such that $q\in D_1\subset\subset D_2\subset\subset D$.  Let 
$V=g(D_2)$.  By Lemma 4.1, $V$ is a locally closed open
Riemann surface.  Since every open Riemann surface is
Stein (cf. [Na2, Thm 3.10.13]) and every holomorphic
line bundle of an open Riemann surface is trivial (cf.
[For, Thm 30.3]), it follows from [Siu, Cor.1] that there exists  
a biholomorphic mapping $\Psi$ from  an open   
neighborhood $W$ of $V$ to an open neighborhood $U$ of $V\times\{0\}$  
in $V\times \C$ such that $\Psi(g(z, w))=(g(z, w), 0)$ for $(z, w)\in D_2$.  
We may assume that $U\subseteq V\times \Delta$.
Let $\pi_1\colon \Delta\to V$ be the universal covering map.  Let  
$\pi(z, w)=(\pi_1(z), w)$.  Let $\Omega=\pi^{-1}(\Psi(W\cap D))$.  Then 
$$ 
\Omega=\{(z', w')\in \Delta\times \C; \quad \pi(z', w')\in U, \ \tilde\rho (z', 
w')< 0 \}
$$ 
where $\tilde\rho (z', w')=\rho(\Psi^{-1}(\pi(z', w'))$.  It is 
easy to see that $\tilde\rho(z', 0)=\frac{\partial\tilde\rho}{\partial
z'}(z', 0)=0$ and $\frac{\partial\tilde\rho}{\partial
w'}(z', 0)\not=0$.  Therefore $b\Omega$ is smooth and hence Levi
flat in a neighborhood of $M=\Delta\times \{0\}$.

Fix a preimage $p' \in \pi^{-1}(\Psi(p))$.  After
a unitary transformation, we may assume that $p'=(0, 0)$. 
Since $\pi$ is locally one-to-one, there exist  unique liftings $q'_j$ of 
$\Psi(q_j)$ for sufficiently large $j$ such that $q'_j\to p'$. 
Since $D_1$ is simply-connected and $g_j(D_1)\subseteq 
W\cap D$ for sufficiently large $j$,  there exist unique liftings
$\tilde g_j$ and $\tilde g 
\colon D_1\to \Omega$ of $\Psi(g_j)$ 
and $\Psi(g)$ respectively such that $\tilde g_j(q)=q'_j$ 
and $\tilde g(q) =p'$.   Choose $\delta\in (0, 1)$ sufficiently 
closed to 1 such that $\tilde g(D_1)\subset\subset \Delta_\delta
\times \{0 \}$.  Let $\eps >0$ and $U_{\delta\eps}
=\Delta_\delta\times\Delta_\eps$.  It follows from Lemma 2.1 that
$$    
\align    
\frac{M^C_{D_1}(q)}{M^E_D(q)}    
&\ge    
\frac{|\det\tilde g'_j (q)|^2 M^C_{\tilde    
g_j(D_1)}(q'_j)}{|\det g'_j (q)|^2 M^E_D(q_j)}     
\ge    
\frac{|\det\tilde g'_j (q)|^2 M^C_{\Omega\cap U_{\delta\eps}}(q'_j)}   
{|\det g'_j (q)|^2 M^E_D(q_j)} \\   
&\ge    
\frac{|\det \tilde g'_j (q)|^2 M^C_{\Omega\cap U_{\delta\eps}}(q'_j)}   
{|\det g'_j (q)|^2 M^E_{\pi(\Omega\cap U_{\delta\eps})}(q_j)}    
=\frac{|\det\tilde g'_j (q)|^2 M^C_{\Omega\cap U_{\delta\eps}}(q'_j)}   
{|\det\tilde g'_j (q)|^2 |\det\pi'(q'_j)|^2 
M^E_{\pi(\Omega\cap U_{\delta\eps})}(q_j)}\\     
 &\ge \frac{M^C_{\Omega\cap U_{\delta\eps}}(q'_j)}   
{M^E_{\Omega\cap U_{\delta\eps}}(q'_j)}     
\endalign    
$$    
By Lemma 4.2, 
$M^C_{D_1}(q)/M^E_D(q) \ge 1$.  
Since we can exhaust $D$ by relatively compact simply-connected 
subdomains (Lemma 4.3), we obtain that 
$M^C_{D}(q)/M^E_D(q) \ge 1$.  It then follows from
Lemma 2.1(4) that $D$ is biholomorphic to the bidisc. \qed

\bigskip

\subhead \S 5.  Proof of the Main Theorem, Part II \endsubhead

In this section, we prove the main theorem when $g(D)\cap\CS\not=
\emptyset$.    

Let $p=g(q)\in g(D)\cap\CS$. 
Then there exists a neighborhood $U$ of 
$p$ such that $U\cap D =\{ (z, w)\in U; \quad \rho_1 < 0, \rho_2 < 0 \}$ where
$\rho_1$ and $\rho_2$ are functions in a defining system of $D$. Let
 $S_1$ and $S_2$ be the corresponding defining hypersurfaces.
It follows from Lemma 3.2 that $g(D)\subseteq S_1\cap S_2$.
Since $\partial\rho_1(p)$ and $\partial\rho_2(p)$
are linearly independent over $\C$, after an affine linear transformation, 
we may assume that $p=(0, 0)$,
$\partial \rho_1 (p)=(0, 1)$, and $\partial \rho_2 (p) =(1, 0)$.

\proclaim{Lemma 5.1} The boundary point $p$ is a local peak point of $D$.
\endproclaim

\demo{Proof}  Assume that the leaves $L_1$  and 
$L_2$ of the foliations of $S_1$ and $S_2$ 
through $p$ are given locally by $\varphi(\zeta)=(\varphi_1(\zeta),
\varphi_2(\zeta))$ and $\psi(\zeta)=(\psi_1(\zeta),
\psi_2(\zeta))$ respectively.  After reparametrization, we
may assume that $\varphi(0)=\psi(0)=p$, $\varphi'(0)=(1, 0)$,
and $\psi'(0)=(0, 1)$.  After a change of coordinates of form
$(z, w)\mapsto (z', w')$ where $z'=z-\psi_1(\psi^{-1}_2(w))$ and
$w'=w -\varphi_2(\varphi^{-1}_1(z))$, we may assume that
$L_1$ and $L_2$ are given locally
by $\Delta_{\delta_1} \times \{ 0 \}$
and $\{ 0 \}\times \Delta_{\delta_2} $ respectively.  
Thus, $S_1$ and $S_2$ are defined by 
$$
\Re w + \Im w \cdot O (|z| + |\Im w|) = 0 \quad\text{and}\quad
\Re z + \Im z \cdot O (|w| + |\Im z|) = 0
$$
respectively.  After possible shrinking of $U$, we have
$$
D\cap U \subseteq \{ \Re w - |\Im w| < 0 \} \times \{ \Re z - |\Im z| < 0 \}.  
$$   
Let $f_p(z, w) = \exp (-(-z)^{2\over 3} - (-w)^{2\over 3})$ where the
cubic root takes the principal branch. 
Then $f_p$ is a local peak function at $p$. \qed
\enddemo

It follows from the maximal principle that $bD$ is variety free at $p$.
Therefore $g(D)=\{p\}$.

\proclaim{Lemma 5.2}Let $U_\eps=\Delta_\eps\times \Delta_\eps$.  Then
 for any sequence $\{ q_j\}$ in $D$ that tends to $p$,
$$ 
\varlimsup_{\eps\to 0}\varlimsup_{j\to \infty}   
\frac{   M^C_{D\cap U_{\eps}}(q_j) }    
{  M^E_{D\cap U_{\epsilon}}(q_j) }  =1     .
$$    
\endproclaim

\demo{Proof} The proof is similar to that of Lemma 4.2 in nature.
We first project $q_j$ to the defining hypersurfaces $S_1$ and
$S_2$.  We then construct a new coordinate system such that 
the leaves through the projections of $q_j$ are
discs on the coordinate planes and the outward normal 
of $bD$ is constant for points on each leaf.  We then compare
the invariant measures of $D\cap U_\eps$ with those  of a product
of two cones.  
We provide details as follows.

As in the proof of Lemma 5.1, we can choose a coordinate system
in a neighborhood of $p$ such that the leaves of the foliations
of $S_1$ and $S_2$  through $p$ are locally given by
$L^1=\Delta_{\eps_0}\times \{ 0\}$ and $L^2=\{0\}\times\Delta_{\eps_0}$
for some $\eps_0 >0$.  By Theorem 3.3,  we can further assume that
the outward normal  of $S_1$ is the positive $\Re w$-axis
for points on $L^1$ and the outward normal direction 
of $S_2$ is the positive $\Re z$-axis
for points on $L^2$.  For $\eps\in (0, \eps_0)$, let 
$\Phi(t, \zeta)=(\zeta, \varphi(t, \zeta))$ be the diffeomorphism
from $(-1, 1)\times \Delta_\eps$ onto a neighborhood of $\Delta_\eps
\times \{ 0 \}$ on $S_1$ as constructed in Lemma 3.1.  Choose $t_j$
such that $t_j\to 0$ and $L^1_j=\Phi(t_j, \Delta_\eps)$ is the leaf
of the foliation of $S_1$ that passes through the projection 
of $q_j$ onto $S_1$ in the direction of the positive $\Re w$-axis.  Let 
$v^1_j(z)=\arg \frac{\partial \rho_1}{\partial w}(z, \varphi(t_j, z))$.
By Theorem 3.3, $v^1_j(z)$ is harmonic on $\Delta_\eps$.  Let $u^1_j(z)$
be its harmonic conjugate such that $u^1_j(0)=0$.  Let $h^1_j(z)
=\exp(-u^1_j(z) + iv^1_j(z))$.  Define $F^1_j\colon (z, w)\to (z', w')$
by $z'=z$ and $w'=(w-\varphi(t_j, z))h^1_j(z)$.   Let 
$\tilde\rho_1(z', w')=\rho_1((F^1_j)^{-1} (z', w'))$ and 
$\tilde\rho_2(z', w')=\rho_2((F^1_j)^{-1} (z', w'))$.  Then
$$
\gather
\tilde\rho_1(z', 0)=0, \ \  \frac{\partial\tilde\rho_1}{\partial z'}(z', 0)=0, \
\ \frac{\partial\tilde\rho_1}{\partial w'}(z', 0)=
|\frac{\partial \rho_1}{\partial w}(z', \varphi(t_j, z'))|e^{u_j(z')} >0 ; \\ 
\tilde\rho_2(0, w')=0, \ \  \frac{\partial\tilde\rho_2}{\partial w'}(0, w')=0, \ 
\
\frac{\partial\tilde\rho_2}{\partial z'}(0, w')=\frac{\partial\rho_2}{\partial 
z}
(0, w'/h^1_j(0) + \varphi(t_j, 0)) > 0, \endgather
$$
for $|z'|< \eps$ and $|w'/h^1_j(0) + \varphi (t_j, 0)| < \eps_0$.
Therefore, in the $(z', w')$-coordinates, $L^1_j=\Delta_\eps\times\{ 0\}$
and the outward normal of $S_1$ is the positive $\Re w'$-axis for points
on $L^1_j$.  Moreover, $L^2=\{(0, w'); \ |w'/h^1_j(0) + \varphi(t_j, 0) | < 
\eps_0 \}$
and the outward normal of $S_2$ is the positive $\Re z'$-axis for points
on $L^2$.   

Since $h^1_j(0)\to 1$ and $\varphi(t_j, 0)\to 0$, $F^1_j(L^2)
\supseteq \{(0, w'); \ \ |w'|\le \eps \}$ and
$\Delta_\eps\times\Delta_{\eps/2}\subseteq
F^1_j(U_\eps) \subseteq  \Delta_\eps\times\Delta_{2\eps}$ for sufficiently
large $j$.  

Now let $\Psi(s, w')$ be the diffeomorphism from $(-1, 1)\times \Delta_\eps$
onto a neighborhood of $\{0\}\times \Delta_\eps$ on $F^1_j(S_2)$ as constructed
in Lemma 3.1.   Choose $s_j$
such that $s_j\to 0$ and $L^2_j=\Psi(s_j, \Delta_\eps)$ is the leaf
of the foliation of $F^1_j(S_2)$ that passes through the projection 
of $q_j$ onto $S_2$ in the direction of the positive $\Re z'$-axis.  Let 
$v^2_j(w')=\arg \frac{\partial \rho_2}{\partial z'}(\psi(s_j, w'), w')$.
Let $u^2_j(w')$
be its harmonic conjugate on $\Delta_\eps$
such that $u^2_j(0)=0$.  Let $h^2_j(w')
=\exp(-u^2_j(w') + iv^2_j(w'))$.  Define $F^2_j\colon (z', w')\to (z'', w'')$
by $z''=(z'- \psi(s_j, w'))h^2_j(w')$ 
and $w''=w'$.   Let 
$\hat\rho_1(z'', w'')=\rho_1((F^2_j)^{-1} (z'', w''))$ and 
$\hat\rho_2(z'', w'')=\rho_2((F^2_j)^{-1} (z'', w''))$.  Then
$$
\gather
\hat\rho_1(z'', 0)=0, \  \frac{\partial\hat\rho_1}{\partial z''}(z'', 0)=0, \
\frac{\partial\hat\rho_1}{\partial w''}(z'', 0)=
\frac{\partial\tilde \rho_1}{\partial w'}
(z''/h^2_j(0) + \psi(s_j, 0), 0) >0 ; \\ 
\hat\rho_2(0, w'')=0, \  \frac{\partial\hat\rho_2}{\partial w''}(0, w'')=0, \
\frac{\partial\hat\rho_2}{\partial z''}(0, w'')=|\frac{\partial\tilde
\rho_2}{\partial z'}
(\psi(s_j, w'), w')|e^{u^2_j(w')} > 0, \endgather
$$
for $|z''/h^2_j(0) + \psi (s_j, 0)| < \eps$
and $|w''|< \eps$.  Therefore, in the $(z'', w'')$-coordinates,
$L^1_j$ and $L^2_j$ are discs on coordinate planes and the
outward normal of $S_1$ (and $S_2$) is the
positive $\Re w''$-axis ($\Re z''$-axis respectively) for 
points on $L^1_j$ ($L^2_j$ respectively).  Let $F_j=F^2_j \circ F^1_j$.
Then for sufficiently large $j$,  $\Delta_{\eps/2}\times\Delta_{\eps/2}
\subseteq F_j(D\cap U_\eps) \subseteq \Delta_{2\eps}\times\Delta_{2\eps}$.
Furthermore, there exists a constant $C>0$ such that 
$$
\gather
\{|z''|< \eps/2; \ \Re z'' + C |\Im z''|^2< 0\}
\times  
\{|w''|< \eps/2; \ \Re w'' + C |\Im w''|^2< 0\}\\
\subseteq F_j(D\cap U_\eps)=\{(z'', w'')\in F_j(U_\eps); \ \
\hat\rho_1(z'', w'') < 0,  \hat\rho_2(z'', w'') < 0 \}
\subseteq\\
\{|z''|< 2\eps; \ \Re z'' - C |\Im z''|^2< 0\}
\times  
\{|w''|< 2\eps; \ \Re w'' - C |\Im w''|^2< 0\} 
\endgather
$$
for sufficiently large $j$.  Thus for any two angles
${\theta_1}\in (0, \pi/2)$ and ${\theta_2}\in (\pi/2, \pi)$, we have
$\Gamma^{\frac{\eps}{2}}_{\theta_1}
\times \Gamma^{\frac{\eps}{2}}_{\theta_1}\subseteq
F_j(\Omega\cap U_{\eps}) \subseteq \Gamma^{2\eps}_{\theta_2}
\times\Gamma^{2\eps}_{\theta_2}$ 
provided $\eps$ is sufficiently small.  Let $q''_j=(z''_j,
w''_j)= F_j(q_j)$. Then
$$
\align
1\ge \frac{   M^C_{\Omega\cap U_\eps}(q_j) }    
{  M^E_{\Omega\cap U_\eps}(q_j) }   
&\ge
\frac{   M^C_{\Gamma^{2\eps}_{{\theta_2}}
\times \Gamma^{2\eps}_{{\theta_2}}}(q''_j) }    
{  M^E_{ \Gamma^{\frac{\eps}{2}}_{{\theta_2}} 
\times \Gamma^{\frac{\eps}{2}}_{{\theta_1}}}(q''_j) } \\    
&= \frac{ M^C_{\Gamma^{2\eps}_{{\theta_2}}}(z''_j) }
{ M^E_{\Gamma^{\frac{\eps}{2}}_{{\theta_1}}}(z''_j) }\cdot  
 \frac{ M^C_{\Gamma^{2\eps}_{{\theta_2}}}(w''_j) }
{ M^E_{\Gamma^{\frac{\eps}{2}}_{{\theta_1}}}(w''_j) } . \\
\endalign
$$
It follows from Lemma 2.2 that the last term can be chosen to be as 
close to 1 as we wish provided
$j\to \infty$, $\eps\to 0$, ${\theta_1}\to (\pi/2)^-$, and 
${\theta_2}\to (\pi/2)^+$, the lemma is now proved. \qed
\enddemo

We now prove the main theorem in this case. We use the notation
and setup at the beginning of the section.  Let $q\in D$ and $q_j=
g_j(q)$.  Let $D_1$ be relatively compact subdomain of $D$ containing
$q$.  For any $\eps >0$, since $g(D)=\{p\}$,
 $g_j(D_1)\subseteq D\cap U_\eps$ for sufficiently large $j$.
Therefore, by Lemma 2.1
$$
\frac{M^C_{D_1}(q)}{M^E_D(q)} =\frac{M^C_{g_j(D_1)}(q_j)}{M^E_D(q_j)}
\ge \frac{M^C_{D\cap U_\eps}(q_j)}{M^E_{D\cap U_\eps}(q_j)}.
$$
By Lemma 5.2, $M^C_{D_1}(q)/M^E_D(q) \ge 1$.  
Exhausting $D$ by $D_1$, we have  
$M^C_{D}(q)/M^E_D(q) \ge 1$.  It then
follows from Lemma 2.1 (4) that $D$ is biholomorphic to
the bidisc.  Note that in this part of the proof 
the simply-connected condition on $D$ is not needed. \qed

\bigskip

%
%
%
%
%
%

\enddocument